# Aspectual Principle, Benford's Law and Russell's Paradox

**Victor P. Novikov**

Nalichnaya 36, 199226
Saint-Petersburg, Russian Federation
vpn@lfp.spb.su

**Abstract**

In the present article Benford's law and Russell's paradox are explained by means of Aspectual Principle.

## 1. Introduction

The different aspects of studying the objects are being traditionally applied for a great variety of sciences. For example, the taxonomy of one and the same biological objects can be based both on morphological and genetic characteristics, that is some aspectual approach is applied to the problems of classification. Another example from physics is that one and the same elementary particle can be treated both as a wave and a corpuscle. In this case the various mathematical tools are used for the investigation of objects within the limits of a single class. Thus aspectual principle (AP) might consist in constructing a new classification or changing the attitudes within the same class of objects. In other words two ways of applying the aspectual principle can be used.

The first variant is based on the extension of the old classification of object. That is something happens in the form of changing the problem as viewed from outside. Therefore such an option can be called the external aspectual principle that is e-AP (external AP). Under the second option the approaches of solving the problems are changing inside the immutable class of objects to be studied. The classification itself does not change (does not widen) but the approach to this class as a whole is changing. The transformation of instrumental methods is taking place, in other words a view to the problem within the system itself is changing. That variant of aspectual principle can be called internal one i.e. i-AP (internal AP).

It should be noted that both types of AP have played a great role in the course of mathematics historical development. Every time when a regular internal crisis had been brewing in mathematics, it was the applying of aspectual principle that made it possible to overcome the contradictions successfully. For example, under the influence of e-AP an extension of classification of numbers was taking place and fractions, irrational and complex numbers have been included into mathematics's arsenal. And the use of i-AP allowed to take a fresh view at those objects of research, i.e., to apply some new instrumental methods. So new kinds of mathematics - algebra, geometry and mathematical logic have come into being.

In such a manner both e-AP and i-AP have historically brought revolutionary changes into mathematics. Every time their application solved some global problems of mathematics development and opened up new mathematical horizons. However e-AP and i-AP have not lost their actuality at present time yet, and they often can assist in solving a number of modern complex mathematical problems. As an example of applying the aspectual principle let us examine one interesting task which is called Benford's law.

## 2. The Problem of Uneven Distribution of First Digits

In 1938 American physicist F. Benford discovered the "Law of Anomalous Numbers". He has analyzed about 20,000 numbers contained in the tables, among which there were some data about the surface area of 335 rivers, as well as specific heat and molecular weight of thousands of chemical compounds and even home numbers of the first 342 persons listed in a biographical directory of American scientists. As a result of that tremendous work he has found a surprising regularity. It turned out that the probability of occurrence **«1»** as the first digit was more than six times greater than that of **«9»**. The probability of occurrence of any other figures was between those extreme values.

**Table 1.** Benford's statistical data

| The first digit of the number | Probability of occurrence as the first digit |
|---|---|
| 1 | 0,306 |
| 2 | 0,185 |
| 3 | 0,124 |
| 4 | 0,094 |
| 5 | 0,080 |
| 6 | 0,064 |
| 7 | 0,051 |
| 8 | 0,049 |
| 9 | 0,047 |
| **Total** | **1,000** |

Since then more than 70 years passed and confirmations of this law were found for a huge multitude of different sets of numbers. For example, the first digits of the magnitude of the population and areas of different countries of the whole world as well as a lot of numeric sequences obey that distribution law. In this way the units also compose about 30% in the following sequence of first digits of powers of two ($2^0$=**1**, $2^1$=**2**, $2^2$=**4**, $2^3$=**8**, $2^4$=**1**6, $2^5$=**3**2, $2^6$=**6**4, $2^7$=**1**28, $2^8$=**2**54, $2^9$=**5**08, $2^{10}$=**1**016, $2^{11}$=**2**032, $2^{12}$=**4**064. . . ).

**1, 2, 4, 8, 1, 3, 6, 1, 2, 5, 1, 2, 4…** (1)
**Sequence of the first digits of powers of two**

In so doing the probability of occurrence of any single-value digit as the first digit, as it is seen from above shown table 1, correlates with **the position** of those digits into the first tens of natural numbers - **1,2,3,4,5,6.7,8,9.** Therefore it is rather interesting to consider such a situation in other positional number systems. For the best illustration, we shall take some extreme cases.

Of course, to the simple binary positional number system the sequence of the first digits of powers of two looks like as a simple infinite sequence of units –

**1, 1, 1, 1, 1, 1, 1….** (2)
**Sequence of the first digits of powers of two in the binary system**

Thus the probability of occurrence of the unit as the first digit of this number system is equal to 100%.

An alternative limiting version of the number system consists of infinite natural numbers, in which each number has got its own unique figure. Of course, it is a hypothetical non-positional system, which does not have any analogy in the history of mathematics. It may be considered as a numberless non-positional system or a degenerate variant of the positional system where there is only one <u>infinite</u> level represented by a natural close. For example the natural numbers can be shown in the following way into this notation system:
**1,2,3,4,5,6,7,8,9, a, b, c, d, e, f, g, h, i, j, k, l, m, n, o, p, q, r, s, t, u, v, w, x, y, z ….∞** (3)

In such a notation system any number of sequences of powers of two (1) is presented by a new unique number; hence the probability of appearance of that number is tending to zero. Of course, it applies to the probability of occurrence of a unit - it appears only once ($2^0=1$), and, as other numbers, it is tending to zero.

It should be noted that in this regard the decimal system is an intermediate variant. Other notations are also quite fit to this sequence. As shown in the Table 2, the probability of appearing the unit as the first digit of powers of two (1) gradually decreases in different number systems. Thus in the ternary notation the probability of appearing the unit as the first digit is equal to about 70%, and in duodecimal one it is about 20%.

**Table 2.** Probability of occurrence of the unit as the first digit in the sequence of powers of two in different number systems

| Base of number system | Probability of occurrence of «1» as the first digit |
|---|---|
| 2 | 1,00 |
| 3 | 0,70 |
| 4 | 0,65 |
| 5 | 0,62 |
| 6 | 0,55 |
| 7 | 0,50 |
| 8 | 0,44 |
| 9 | 0,38 |
| 10 | 0,31 |
| 11 | 0,25 |
| 12 | 0,20 |
| ∞ | 0,00 |

Thus, the use of aspectual principle can explain the appearance of "The Law of Uneven Distribution of First Digits". Some new regulations have been displayed when using some different number systems within the limits of natural figures and it allows us to conclude that Benford's law is not determined by any intrinsic properties of natural numbers, but by an arbitrary choice of the number system only. With the same probability a ternary or duodecimal and even sexagesimal number system, but not decimal one, could be applied to modern mathematics. And then the Benford's law of distribution of first digits would look very differently in each separate case. But it should be noted that under any system of notation the distribution of possibilities of first digits depends, to some extent, on their position to be placed in the first digit natural numbers. Application of the i-AP at the present example gave an opportunity to look at that problem from an absolute new angle of view. Or you can express it in another way - the use of i-AP is similarly to mental transfer of a research object to another class or creation of a parallel classification cell.

However, in more complex cases, a simple transfer of the object to the adjacent parallel cell does not have any effect. Some stronger and more radical variant of aspectual principle requires to be applied that is e-AP. Because it is necessary to look at the problem from the side, in other words, to move beyond the initial limits of classification to a higher generalized level or, in extreme cases, to create a new classification, where the object of research will be included as a special case. As an example of using e-AP a well-known Russell's paradox (also known as Russell's antinomy) may be examined.

# 3. Russell's Paradox

The solution to this antinomy fully demonstrates the beauty and elegance of applying the e-AP because it is a combination of clearness and inevitability. The paradox itself has got a long prehistory with the details to be widely represented in various publications. At the very beginning of the twentieth century Bertrand Russell, in a letter to Frege drew his attention to the incorrect use of the term "class of all classes" to be applied in the set theory, which lies at the foundation of arithmetic. Russell's antinomy is formulated as follows:
"Let R be the set of all sets that are not members of themselves. If R qualifies as a member of itself, it would contradict its own definition as a set containing sets that are not members of themselves. On the other hand, if such a set is not a member of itself, it would qualify as a member of itself by the same definition."
Symbolically:

$$\text{let } R = \{x \mid x \notin x\}, \text{ then } R \in R \iff R \notin R$$

There Russell illustrates an example of this paradox, which he called as "paradox of a barber." Once in a small town a barber declared that he would shave everyone who did not shave himself. In so doing he boasted that he has no equal as for hairdressing but once he thought whether he should shave himself. On the one hand he should not do this because he shaves only others. But if he does not shave himself, then he will come in the number of those who do not shave themselves, and therefore he gets a possibility to shave himself.

Trying to understand the reason of the paradox Russell and Whitehead wrote in their mutual book "Principia mathematica" (1910-1913): "The analysis of paradoxes determines that all of them have got a vicious circle as a resource. That vicious circle arises from the adoption that a great number of objects are able to contain an element which might be only determined by the set as a whole. "In that case, returning to the example of the hairdresser, we can say that only two assertions in separate would be non-contradictory but not one assertion in common. Or this barber would be able to shave all clients except himself, or he with the rest of the residents together is a subject to be shaved by another barber. That is there is no possibility of merging both sets into a single class. Strictly speaking from the standpoint of the aspectual principle, this is the only right solution of the paradox that is to create one class instead of two ones. A class of clients and the class of practicing barbers can exist only in separate condition, who may be customers in one class, and the barbers in the other one.

Thus the extension of the classification with the help of e-AP completely solves the problem by creating a new, more generalized version of the classification where the object of study is included as a special case. Contradiction disappears simply by definition, since it had been created artificially (in the limits of starting classification) and it was also solved artificially (in a new classification).

Apparently, in most cases, when in fact there is an unavoidable contradiction, e-AP can and should be applied. In mathematics almost all of the contradictions have been caused just by unsuccessful or simply improper classification, and then the only way is to create a new classification. Of course, you must use it carefully and don't propagate the essence beyond measure! However you needn't hold to the idea that the whole mathematical classification was anticipated and sent down to us by the ancient Egyptian god of calculation namely Thoth, so it is said to be sacred by definition. But even if you accept it, now it is not clear where the messages of Thoth are ending and fabrications of our contemporaries begin. Besides, it is hardly possible to believe that a single immutable final classification of mathematics is able to exist in general which would be suitable for all occasions. These doubts could be solved only by the god Thoth himself, having been created at the times of the powerful pharaohs.

# 4. Conclusion. Aspectual principle in gnosiology (epistemology)

Above presented thoughts allow to draw a conclusion that the aspectual principle in mathematics has not lost its value nowadays. Its application makes possible solving the tasks which require either an extension of classification in the direction of synthesis (e-AP) or some transformation of mathematical tools in the limits of parallel class (i-AP). In so doing it should be noted that the use of the aspectual principle allows to give a fresh look and evaluate the part of mathematics in a unique science building and the entire learning process. Practically, in most cases, i-AP can be applied to solve some mathematical tasks. By the way, as an object of investigation remains unchanged here, such a procedure can be called the operation of **<u>aspectual symmetry</u>**.

Let's take a notion of a real number as one easy example. Concerning that subject H.Weyl (1885-1955) wrote: "The system of real numbers is similar to Janus: on one side it is a set of algebraic operations as "+" and "-" and their inverses and on the other one it is a continual variety, their parts are continuously related to each other. The first face of the numbers is algebraic , the second one is topological ". It should be also noted that one and the same object of study can be investigated in biology, chemistry, physics and mathematics. Consequently i-AP is acting in epistemology in that case too.

However, in general, you can learn the world around us not only with the help of instruments, but also, for example, using literature, poetry, music as an object of cognition. In this case, using the e-AP it is going beyond the scope of usual scientific paradigm. You can even spend some parallels between arts and science. For example, Biology (and Geology) can be correlated with the Documentary, Physics - with Romance, and Math with Science Fiction. The last of the comparisons is particularly interesting, since in both cases, the specialists are constructing some abstract and often obscure designs. Thus the aspectual principle is an extremely powerful means for the study of our surroundings.

## References

1. Benford, F. "The Law of Anomalous Numbers." Proc. Amer. Phil. Soc. 78, 551-572, 1938
2. Russell, B. "The principles of mathematics." Cambridge: University press, 1903
3. Russell, B. and Whitehead, A "Principia Mathematica." in 3 vols. Cambridge University Press. 1910–13
4. Weyl, H. "The Classical Groups. Their Invariants and Representations." Princeton University Press. 1939